\documentclass{commat}

\usepackage[all]{xy}

\newcommand{\Z}{\mathbb{Z}}
\DeclareMathOperator{\Id}{Id}
\DeclareMathOperator{\SL}{SL}
\DeclareMathOperator{\PSL}{PSL}

\title{Combinatorial description of the principal congruence subgroup \(\Gamma(2)\) in \(\SL(2,\Z)\)}

\author{Flavien Mabilat}

\affiliation{
\address{%
Laboratoire de Mathématiques de Reims, 
UMR9008 CNRS et Université de Reims Champagne-Ardenne, 
U.F.R. Sciences Exactes et Naturelles Moulin de la Housse - BP 1039 51687 Reims cedex 2, 
France}
\email{flavien.mabilat@univ-reims.fr}
}

\abstract{
We characterize sequences of positive integers
\((c_{1},c_{2},\dots,c_{n})\)
for which the \(2\times2\) matrix \(
\begin{pmatrix}
   c_{1} & -1 \\
    1 & 0
   \end{pmatrix}
\cdots
   \begin{pmatrix}
   c_{n} & -1 \\
    1 & 0
    \end{pmatrix}
\) belongs to the principal congruence subgroup of level 2 in \(\SL(2,\Z)\).
The answer is given in terms of dissections of a~ convex \(n\)-gon into a~ mixture
of triangles and quadrilaterals.
}

\keywords{modular group, polygon, quiddity.}

\msc{05A05}

\DOI{10.46298/cm.9039}

\VOLUME{30}

\firstpage{13} 
\lastpage{23}

\begin{paper}

\section{Introduction}

The classical \textit{modular group}
\[
\SL(2,\Z) :=
\left\{
\begin{pmatrix}
a & b \\
c & d
   \end{pmatrix}
 \;\vert\;a,b,c,d\in\Z,\;
 ad-bc = 1
\right\}
\]
and its quotient by the center, \(\PSL(2,\Z) :=\SL(2,\Z)/\{\pm\Id\}\), play a~ central role in several classical areas, such as the theory of continued fractions, hyperbolic geometry, and the theory of modular forms.  The group \(\SL(2,\Z)\) naturally acts on the upper half-plane, and perhaps the most remarkable fact about it is that the quotient by this action is the moduli space of elliptic curves (this fact explains the name ``modular group'' due to Klein).  The structure of the modular group and its subgroups has been thoroughly studied, see~\cite{R}.  An important class of subgroups, \(\Gamma(N)\), are called \textit{principal congruence subgroups of level N}.  These are defined as follows:
\[
\Gamma(N) :=
\left\{
A \in \SL(2,\Z)\;\vert\;
 A = \Id\;(\!\!\!\!\mod N)
\right\},
\]
where \(N\) is a~ positive integer.

This note is about a~ relation of the modular group to combinatorics.
The idea is based on the fact that every element of \(\SL(2,\Z)\)
has a~ (canonical) presentation (i.e. a~ description by means of generators and relations) by a~ sequence of \textit{positive integers}.
This has been known for a~ long time (cf.~\cite{Zag}), but
started to be exploited only very recently; see~\cite{O},\cite{MO}.
One uses a~ general principle that positive integers
must count some (geometric, combinatorial, etc.) objects.

Our approach is based on the work of Coxeter~\cite{Cox} and Conway-Coxeter~\cite{CoCo},~\cite{CoCo2}. Coxeter and Conway used the notion of \textit{frieze pattern} to characterize sequences of positive integers
\((c_1,\ldots,c_n)\) such that
\begin{equation}\label{CoCEqI}
\begin{pmatrix}
   c_{1} & -1 \\
    1 & 0
   \end{pmatrix}
\begin{pmatrix}
   c_{2} & -1 \\
    1 & 0
   \end{pmatrix}
   \cdots
   \begin{pmatrix}
   c_{n} & -1 \\
    1 & 0
    \end{pmatrix}
= -\Id,
\end{equation}
and that satisfy an extra condition of \textit{total positivity}, formulated as the positivity of the entries of
Coxeter's \textit{frieze pattern}, or just \textit{frieze} for short.
All positive solutions of equation~\eqref{CoCEqI} were classified in~\cite{O}.
For a~ detailed explanation of the total positivity property, see~\cite{O},\cite{MO} (and also~\cite{CO}).

Our goal is twofold.
We give a~ short overview of the combinatorial approach to the modular group,
that we believe should be better known.
We prove a~ new theorem that gives a~ combinatorial description
of \(\Gamma(2)\), the principal congruence subgroup of level 2.

\section{Sequences of positive integers}

The group \(\SL(2,\Z)\) has two generators whose standard choice is
\[
  T :=
\begin{pmatrix}
   1 & 1 \\
    0 & 1
   \end{pmatrix},
\qquad
 S :=
\begin{pmatrix}
   0 & -1 \\
    1 & 0
   \end{pmatrix}.
 \]
These generators satisfy the relations: \(S^{2} = {(TS)}^{3} = -\Id\),
and this is a~ complete set of relations in \(\SL(2,\Z)\).
This classical fact can be found in many textbooks; for a~ particularly elementary proof, see~\cite{Alp}.
It readily implies that every element \(A\) of \(\SL(2,\Z)\) can be written,
for some positive integer \(n\), in the form
\begin{equation}\label{DecMEq}
A = \pm T^{c_{1}}S\,T^{c_{2}}S\cdots T^{c_{n}}S = \pm
\begin{pmatrix}
   c_{1} & -1 \\
    1 & 0
   \end{pmatrix}
\begin{pmatrix}
   c_{2} & -1 \\
    1 & 0
   \end{pmatrix}
   \cdots
   \begin{pmatrix}
   c_{n} & -1 \\
    1 & 0
    \end{pmatrix}
,
\end{equation}
where \(c_{1},\ldots,c_{n}\) are positive integers; see~\cite{R},\cite{Zag},\cite{O}, and the explanation below.
We will use the notation
\[
M(c_1,\ldots,c_n) :=
\begin{pmatrix}
   c_{1} & -1 \\
    1 & 0
   \end{pmatrix}
\begin{pmatrix}
   c_{2} & -1 \\
    1 & 0
   \end{pmatrix}
   \cdots
   \begin{pmatrix}
   c_{n} & -1 \\
    1 & 0
    \end{pmatrix}
,
\]
For the generators one easily checks
\[
T = -M(2,1,1),
\qquad
T^{-1} = -M(1,1,2,1),
\qquad\hbox{and}\qquad
S = -M(1,1,2,1,1).
\]
A decomposition \(A = \pm M(c_1,\ldots,c_n)\) with each \(c_i\) positive
can then be obtained for every chosen \(A\) by concatenation
of any expression of \(A\) in terms of the generators.

The decomposition \(A = \pm M(c_1,\ldots,c_n)\) is obviously not unique
(even though a~ canonical, shortest expression was suggested in~\cite{MO}).
The first natural problem is thus to consider the equation
\begin{equation}\label{TheEqn}
M(c_1,\ldots,c_n) = \pm\Id,
\end{equation}
and look for a~ combinatorial description of the sequences of positive integer solutions.
In other words, this problem is to give an explicit combinatorial description of
relations in \(\PSL(2,\Z)\).
This problem was studied in~\cite{CoCo},~\cite{CoCo2},~\cite{O}; see also~\cite{ARS},~\cite{MO},~\cite{Hen}
and Section~\ref{PDSec} below.
It turns out that equation~\eqref{TheEqn} is related to triangulations of \(n\)-gons and to more
sophisticated polygon dissections.

\section{The main result of this paper}

We will generalize equation~\eqref{TheEqn} and describe the sequences of positive integers
\((c_1,\ldots,c_n)\) for which
\begin{equation}\label{ThePbl}
M(c_1,\ldots,c_n) \in \Gamma(2),
\end{equation}
where \(\Gamma(2)\) is the principal congruence subgroup of \(\SL(2,\Z)\) of level \(2\), see the introduction.

Similarly to the case of equation~\eqref{TheEqn}, the property of being a~ solution of equation~\eqref{ThePbl} is cyclically invariant
(i.e., if an \(n\)-tuple \(({c_{1}},\ldots,{c_{n}})\) is a~ solution of equation~\eqref{ThePbl},
then \(({c_{n}},{c_{1}},\ldots,{c_{n-1}})\) is also a~ solution).
It is thus often convenient to consider, instead of an \(n\)-tuple
\(({c_{1}},\ldots,{c_{n}})\), an \(n\)-periodic infinite sequence \(\left({c_{i}}\right)_{i\in\Z}\).

The solutions of equation~\eqref{ThePbl}
can be formulated in terms of polygon dissections.

\begin{definition}
Let a~ \textit{\((3|4)\)-dissection} be any partition of a~ convex \(n\)-gon into sub-polygons by
pairwise non-crossing diagonals, such that every subpolygon is a~ triangle or a~ quadrilateral.
\end{definition}

\begin{example}\label{TheEx}
Let us give here simple examples:
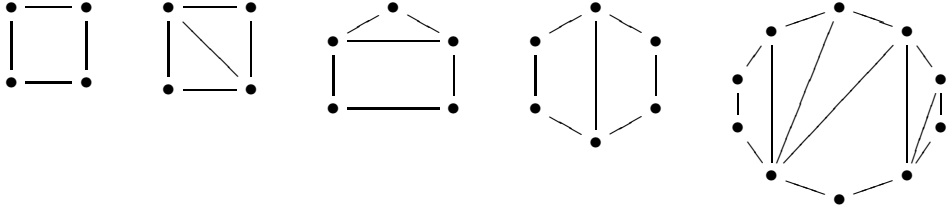
\[
\xymatrix @!0 @R = 0.5cm @C = 0.5cm
{
\bullet\ar@{-}[dd]\ar@{-}[rr]&&\bullet\ar@{-}[dd]
\\
\\
\bullet\ar@{-}[rr]&& \bullet
}
\quad
\xymatrix @!0 @R = 0.55cm @C = 0.55cm
{
\bullet\ar@{-}[dd]\ar@{-}[rr]\ar@{-}[dd]\ar@{-}[rrdd]&&\bullet\ar@{-}[dd]
\\
\\
\bullet\ar@{-}[rr]&& \bullet
}
\quad
\xymatrix @!0 @R = 0.45cm @C = 0.8cm
{
&\bullet\ar@{-}[ld]\ar@{-}[rd]&
\\
\bullet\ar@{-}[dd]\ar@{-}[rr]&&\bullet\ar@{-}[dd]
\\
\\
\bullet\ar@{-}[rr]&& \bullet
}
\quad
\xymatrix @!0 @R = 0.45cm @C = 0.8cm
{
&\bullet\ar@{-}[ld]\ar@{-}[rd]\ar@{-}[dddd]&
\\
\bullet\ar@{-}[dd]&&\bullet\ar@{-}[dd]
\\
\\
\bullet&& \bullet
\\
&\bullet\ar@{-}[lu]\ar@{-}[ru]&
}
\quad
\xymatrix @!0 @R = 0.32cm @C = 0.45cm
 {
&&&\bullet\ar@{-}[llddddddd]\ar@{-}[lld]\ar@{-}[rrd]&
\\
&\bullet\ar@{-}[ldd]\ar@{-}[dddddd]&&&& \bullet\ar@{-}[rdd]\ar@{-}[lllldddddd]\ar@{-}[dddddd]\\
\\
\bullet\ar@{-}[dd]&&&&&&\bullet\ar@{-}[dd]\ar@{-}[ldddd]\\
\\
\bullet\ar@{-}[rdd]&&&&&&\bullet\ar@{-}[ldd]\\
\\
&\bullet\ar@{-}[rrd]&&&& \bullet\ar@{-}[lld]\\
&&&\bullet&
}
\]
\end{example}

For any integer \(a\), we denote by \(\overline{a} :=a+2\mathbb{Z}\) the image of \(a\) under the projection \(\Z\to\Z/2\Z\).  If \(a\) is odd, then \(\overline{a} = \overline{1}\); if \(a\) is even, then \(\overline{a} = \overline{0}\).
The following notion is analogous to that of~\cite{CoCo} and~\cite{CoCo2}.

\begin{definition}
The \textit{quiddity} of a~ \((3|4)\)-dissection of a~ convex \(n\)-gon is the (cyclically ordered) sequence \((\overline{c_1},\ldots,\overline{c_n})\)
of elements of \(\Z/2\Z\), such that for every \(i\)
\[
\overline{c_{i}} =
\left\{
\begin{array}
{ll}\overline{1}, & \hbox{if the number of triangles
adjacent to the \(i\)th vertex is odd};\\
[2pt]
\overline{0}, & \hbox{if the number of triangles
adjacent to the \(i\)th vertex is even}.
\end{array}
\right.
\]
\end{definition}

\begin{example}\label{TheExBis}
The quiddities \((\overline{c_1},\ldots,\overline{c_n})\) of the \((3|4)\)-dissections of Example~\ref{TheEx}
are as follows:
\begin{enumerate}
\item[(a)]
For the first two examples,
\[
(\overline{c_1},\overline{c_2},\overline{c_3},\overline{c_4}) 
= \left(\overline{0},\overline{0},\overline{0},\overline{0}\right)
\quad \textup{and} \quad 
(\overline{c_1},\overline{c_2},\overline{c_3},\overline{c_4}) =
(\overline{0},\overline{1},\overline{0},\overline{1}).
\]
These are the only quiddities for \(n = 4\).
\item[(b)]
For the remaining examples one has, respectively,
\begin{gather*}
(\overline{c_1},\ldots,\overline{c_5})
= (\overline{1},\overline{1},\overline{1},\overline{0},\overline{0}), \quad
(\overline{c_1},\ldots,\overline{c_6})
= (\overline{0},\overline{0},\overline{0},\overline{0},\overline{0},\overline{0}) \\
\textup{and} \quad
(\overline{c_1},\ldots,\overline{c_{10}})
= (\overline{0},\overline{1},\overline{0},\overline{0},
\overline{0},\overline{0},\overline{0},\overline{1},\overline{0},\overline{0}).
\end{gather*}
\end{enumerate}
\end{example}

The following statement is our main result.  It gives a~ combinatorial characterization of the solutions of equation~\eqref{ThePbl} for \(n\geq3\).  Note that the product of elements of \(\SL(2,\Z)\) commutes with the projection
of the entries of matrices to \(\Z/2\Z\).  This allows one to make all the computations in \(\SL(2,\Z/2\Z)\).

\begin{theorem}\label{MaThm}
\mbox{}
\begin{itemize}
\item[(i)] Every quiddity of a~ \((3|4)\)-dissection of an \(n\)-gon is a~ solution of equation~\eqref{ThePbl}.
\item[(ii)] Every solution of equation~\eqref{ThePbl} with \(n\geq3\) is a~ quiddity of a
\((3|4)\)-dissection of an \(n\)-gon.
\end{itemize}
\end{theorem}

This statement is proved in Section~\ref{ProofSec}.

Let us mention that the number of solutions of equation~\eqref{ThePbl}, for a~ fixed \(n\), can be deduced from the main result of~\cite{Sop1} and is given by the Jacobsthal sequence (A001045 in OEIS~\cite{OEIS}).

\section{Relations in \(\PSL(2,\Z)\) and polygon dissections}\label{PDSec}

We give a~ brief overview of the theorems of Conway and Coxeter~\cite{CoCo},~\cite{CoCo2} (see also~\cite{ARS},\cite{Hen})
and Ovsienko~\cite{O}.
The first one relates equation~\eqref{TheEqn} to one of the most classical notion of combinatorics,
namely that of triangulation of an \(n\)-gon, while the second describes all positive integer solutions
of equation~\eqref{TheEqn} in terms of polygon dissections.
This overview will allow us to compare equation~\eqref{TheEqn} and equation~\eqref{ThePbl}.
It also makes the presentation complete.

\subsection{Triangulations and friezes}\label{CoCoSec}

Fix a~ triangulation of a~ convex \(n\)-gon.
Following~\cite{CoCo} and~\cite{CoCo2}, we call a~ \textit{quiddity} of the triangulation the sequence of positive integers \((c_1,\ldots,c_n)\), where \(c_i\) is equal to the number of triangles adjacent to the \(i\)th vertex of the \(n\)-gon.

The theorem of Conway and Coxeter can be formulated in the following way (cf.~\cite[Corollary 2.3]{O}).

\begin{theorem}[\cite{CoCo} and~\cite{CoCo2}]\label{CoCoThm}
\mbox{}

\begin{itemize}
\item[(i)] The quiddity of a~ triangulation of an \(n\)-gon is a~ solution of the equation \(M(c_1,\ldots,c_n) = -\Id\).
\item[(ii)] Every solution \((c_1,\ldots,c_n)\) of equation~\eqref{TheEqn} satisfying the condition
\begin{equation}\label{TheTPEqn}
c_1+c_2+\cdots+c_n = 3n-6
\end{equation}
is the quiddity of a~ triangulation of an \(n\)-gon.
\end{itemize}
\end{theorem}

The simplest examples, with \(n = 3,4,5\), are
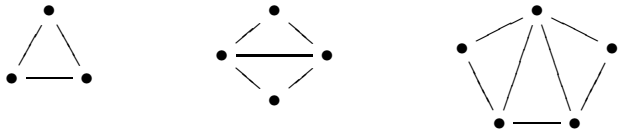
\[
\xymatrix @!0 @R = 0.45cm @C = 0.5cm
{
&\bullet\ar@{-}[ldd]\ar@{-}[rdd]&
\\
\\
\bullet\ar@{-}[rr]&& \bullet
}
\qquad
\qquad
\xymatrix @!0 @R = 0.60cm @C = 0.7cm
{
&\bullet\ar@{-}[ld]\ar@{-}[rd]&
\\
\bullet\ar@{-}[rr]&& \bullet
\\
&\bullet\ar@{-}[lu]\ar@{-}[ru]&
}
\qquad
\qquad
\xymatrix @!0 @R = 0.50cm @C = 0.5cm
{
&&\bullet\ar@{-}[rrd]\ar@{-}[lld]\ar@{-}[lddd]\ar@{-}[rddd]&
\\
\bullet\ar@{-}[rdd]&&&& \bullet\ar@{-}[ldd]\\
\\
&\bullet\ar@{-}[rr]&& \bullet
}
\]
It is easy to see that their corresponding quiddities, namely \((c_1,c_2,c_3) = (1,1,1)\), \((c_1,c_2,c_3,c_4) = (1,2,1,2)\) and \((c_1, c_2,c_3, c_4,c_5) = (1, 3, 1, 2, 2)\), are, indeed, solutions of the equation \(M(c_1,\ldots,c_n) = -\Id\).

The original formulation of Theorem~\ref{CoCoThm} uses the beautiful notion of Coxeter's \textit{frieze}.
Let us recall that a~ frieze is an array of \((n-1)\) infinite rows of positive integers with
the rows \(1\) and \(n-1\) consisting in \(1\)'s.
Every elementary \(2\times2\) ``diamond''
\[
\begin{array}
{ccc}
&b&\\
a&&d\\
&c&
\end{array}
,
\]
of the frieze must satisfy the unimodular rule \(ad-bc = 1\).
Coxeter proved in~\cite{Cox} that the row \(2\) (and \(n-2\)) is an \(n\)-periodic sequence
satisfying equation~\eqref{CoCEqI}.
The Conway-Coxeter theorem of~\cite{CoCo} and~\cite{CoCo2} identifies Coxeter's friezes with triangulations.

Let us give here an example of a~ frieze for \(n = 5\):
\[
 \begin{array}
{cccccccccccccccc}
\cdots&&1&& 1&&1&&1&&1
 \\
[2pt]
&1&&3&&1&&2&&2&&\cdots
 \\
[2pt]
\cdots&&2&&2&&1&&3&&1
 \\
[2pt]
&1&& 1&&1&&1&&1&&\cdots
\end{array}
\]
The \(5\)-periodic sequence \((1,3,1,2,2)\) is the quiddity of a~ triangulation of a~ pentagon.

For a~ survey on friezes, see~\cite{Sop}.  Variants of friezes involving other types of polygon dissections can be found in~\cite{BHJ}, \cite{HJ}, \cite{CH}.  Links of friezes to presentations of \(\SL(2,\Z)\) also appeared in~\cite{BR}.
\begin{remark}
Let us mention that~\eqref{TheTPEqn} turns out to be equivalent
to the condition of total positivity, see~\cite[Corollary 2.3]{O};
and it can be formulated in more standard terms of continued fractions and total positive \((2\times2)\)-matrices, see~\cite{MO}.
In terms used by Coxeter, this total positivity means that every entry of the frieze is positive.
\end{remark}

\subsection{Complete solution of equation~\eqref{TheEqn}}\label{3dSec}

For \(n\geq6\), there exist many solutions of equation~\eqref{TheEqn}
that cannot be obtained from triangulations of \(n\)-gons.
The complete solution of equation~\eqref{TheEqn} was given in~\cite{O} and led
to the following notion of ``\(3d\)-dissection''.

\begin{definition}\ 
\begin{enumerate}
\item[(i)] A \(3d\)-\textit{dissection} is a~ partition of a~ convex \(n\)-gon
into sub-polygons by means of pairwise non-crossing diagonals,
such that the number of vertices of every sub-polygon is a~ multiple of \(3\).

\item[(ii)] The \textit{quiddity} of a~ \(3d\)-dissection of an \(n\)-gon is the
(cyclically ordered) \(n\)-tuple of positive integers
\((c_1,\ldots,c_n)\) such that \(c_i\) is the number of sub-polygons adjacent to the \(i\)th vertex
of the \(n\)-gon.
\end{enumerate}
\end{definition}

\begin{theorem}[\cite{O}, Theorem 1]\label{SecondMainThm}
Every quiddity of a~ \(3d\)-dissection of an \(n\)-gon is a~ solution of equation~\eqref{TheEqn}.
Conversely, every solution of equation~\eqref{TheEqn} is a~ quiddity of
a \(3d\)-dissection of an \(n\)-gon.
\end{theorem}

The simplest examples of \(3d\)-dissections that are not triangulations are
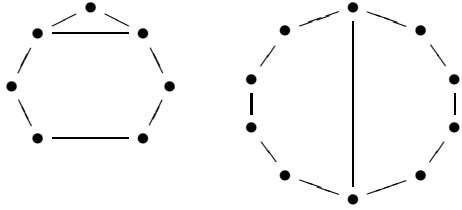
\[
\xymatrix @!0 @R = 0.35cm @C = 0.35cm
{
&&&\bullet\ar@{-}[rrd]\ar@{-}[lld]
\\
&\bullet\ar@{-}[ldd]\ar@{-}[rrrr]&&&& \bullet\ar@{-}[rdd]&\\
\\
\bullet\ar@{-}[rdd]&&&&&& \bullet\ar@{-}[ldd]\\
\\
&\bullet\ar@{-}[rrrr]&&&&\bullet
}
\qquad
\xymatrix @!0 @R = 0.32cm @C = 0.45cm
 {
&&&\bullet\ar@{-}[dddddddd]\ar@{-}[lld]\ar@{-}[rrd]&
\\
&\bullet\ar@{-}[ldd]&&&& \bullet\ar@{-}[rdd]\\
\\
\bullet\ar@{-}[dd]&&&&&&\bullet\ar@{-}[dd]\\
\\
\bullet\ar@{-}[rdd]&&&&&&\bullet\ar@{-}[ldd]\\
\\
&\bullet\ar@{-}[rrd]&&&& \bullet\ar@{-}[lld]\\
&&&\bullet&
}
\]
and the corresponding quiddities are: \((2,1,2,1,1,1,1)\) and \((2,1,1,1,1,2,1,1,1,1)\).

\subsection{Idea of the proof}\label{IdPrSec}

The proof of Theorem~\ref{SecondMainThm} is inductive.
The main idea uses the following ``local surgery'' operations:
\begin{enumerate}
\item[(\(\alpha\))]
\((c_1,\ldots,c_i,c_{i+1},\ldots,c_n)
\mapsto
(c_1,\ldots,c_i+1,\,1,\,c_{i+1}+1,\ldots,c_n)\),

\item[(\(\beta\))]
\((c_1,\ldots,c_{i-1},c_i,c_{i+1}\ldots,c_n)
\mapsto
(c_1,\ldots,c_{i-1},c'_i,\,1,\,\,1,\,c''_i,c_{i+1},\ldots,c_n)\),
\end{enumerate}
where \(c'_i+c''_i = c_i+1\).
One readily checks that the matrix \(M(c_1,\ldots,c_n)\) is invariant under the
operations of type (\(\alpha\)) and changes the sign under the operations of type~(\(\beta\)).

A simple lemma then states that a~ sequence of positive integers \((c_1,\ldots,c_n)\)
satisfying equation~\eqref{TheEqn} always has some entries \(c_i = 1\); cf.~\cite{CoCo},~\cite{CoCo2},~\cite{O}.
This allows one to construct any solution of equation~\eqref{TheEqn} from the
simplest solution \((1,1,1)\).

The inductive step in the proof is based on the observation that the above surgery operations have a~ combinatorial meaning.  Given a~ dissection of an \(n\)-gon, the operation (\(\alpha\)) consists in gluing an additional triangle on the edge \((i,i+1)\), while the operation (\(\beta\)) changes a~ \(3k\)-gon adjacent to the \(i\)th vertex in the dissection into a~ \((3k+3)\)-gon; see~\cite{O}.

\section{Proof of Theorem~\ref{MaThm}}\label{ProofSec}

Our proof of Theorem~\ref{MaThm} is quite similar to that of Theorem~\ref{SecondMainThm}.
We give an inductive procedure of construction of all the solutions of equation~\eqref{ThePbl}.

\subsection{Local surgery}
Consider the following two families of ``local surgery'' operations
 for sequences of elements of \(\Z/2\Z\):

\begin{enumerate}
\item[(a)]
Operations of the first family insert \(\overline{1}\) into the sequence
\((\overline{c_{1}},\overline{c_{2}},\ldots,\overline{c_{n}})\):
\[
(\overline{c_{1}},\ldots,\overline{c_{i}},\overline{c_{i+1}},\ldots,\overline{c_{n}})
\mapsto
(\overline{c_{1}},\ldots,\overline{c_{i}}+\overline{1},\,\overline{1},\,
\overline{c_{i+1}}+\overline{1},\ldots,\overline{c_{n}}).
\]

\item[(b)]
Operations from the second family insert two copies of \(\overline{0}\) between two consecutive positions:
\[
(\overline{c_{1}},\ldots,\overline{c_{i}},\overline{c_{i+1}},\ldots,\overline{c_{n}})
\mapsto
(\overline{c_{1}},\ldots,\overline{c_{i}},\overline{0},\overline{0},\overline{c_{i+1}},\ldots,\overline{c_{n}}).
\]
\end{enumerate}

\noindent
Within the cyclic ordering, the operations (a) and (b) are defined for all
\(1\leq{}i\leq{}n\).  Every operation (a) transforms a~ sequence of \(n\) elements of \(\Z/2\Z\) into a~ sequence of \(n+1\) elements of \(\Z/2\Z\), and every operation (b) transforms a~ sequence of \(n\) elements of \(\Z/2\Z\) into a~ sequence of \(n+2\) elements of \(\Z/2\Z\).

The following statement means that equation~\eqref{ThePbl}
is invariant under the operations (a) and (b).

\begin{lemma}\label{FirstLem}
One has, in the group \(\SL(2,\Z/2\Z)\),
\[
\begin{array}
{rcl}
M(\overline{c_1},\ldots,\overline{c_{n+1}}) &=&
M(\overline{c_{1}},\ldots,\overline{c_{i}}+\overline{1},\,\overline{1},\,
\overline{c_{i+1}}+\overline{1},\ldots,\overline{c_{n}}),\\
[4pt]
M(\overline{c_1},\ldots,\overline{c_{n+1}}) &=&
-M(\overline{c_{1}},\ldots,\overline{c_{i}},\overline{0},\overline{0},\,
\overline{c_{i+1}},\ldots,\overline{c_{n}}).
\end{array}
\]
\end{lemma}

\begin{proof}
An operation of type (a) replaces the matrix \(
\begin{pmatrix}
    \overline{c_{i}} & \overline{1} \\
    \overline{1} & \overline{0}
    \end{pmatrix}
\begin{pmatrix}
    \overline{c_{i+1}} & \overline{1} \\
    \overline{1} & \overline{0}
    \end{pmatrix}
\) by
\begin{align*}
\begin{pmatrix}
    \overline{c_{i}+1} & \overline{1} \\
    \overline{1} & \overline{0}
    \end{pmatrix}
\begin{pmatrix}
    \overline{1} & \overline{1} \\
    \overline{1} & \overline{0}
    \end{pmatrix}
\begin{pmatrix}
    \overline{c_{i+1}+1} & \overline{1} \\
    \overline{1} & \overline{0}
    \end{pmatrix}
&=
\begin{pmatrix}
    \overline{c_{i}} & \overline{c_{i}+1} \\
    \overline{1} & \overline{1}
    \end{pmatrix}
\begin{pmatrix}
    \overline{c_{i+1}+1} & \overline{1} \\
    \overline{1} & \overline{0}
    \end{pmatrix}
\\
&=
\begin{pmatrix}
    \overline{c_{i}c_{i+1}+1} & \overline{c_{i}} \\
    \overline{c_{i+1}} & \overline{1}
    \end{pmatrix}
\\
&=
\begin{pmatrix}
    \overline{c_{i}} & \overline{1} \\
    \overline{1} & \overline{0}
    \end{pmatrix}
\begin{pmatrix}
    \overline{c_{i+1}} & \overline{1} \\
    \overline{1} & \overline{0}
    \end{pmatrix}.
\end{align*}
Therefore, \(M(\overline{c_1},\ldots,\overline{c_{n+1}}) =
M(\overline{c_{1}},\ldots,\overline{c_{i}}+\overline{1},\,\overline{1},\,
\overline{c_{i+1}}+\overline{1},\ldots,\overline{c_{n}})\),
as an element of \(\SL(2,\Z/2\Z)\). An operation of type (b) adds \(-\Id\) in the product that defines \(M_{n}(\overline{c_{1}},\ldots,\overline{c_{n}})\).
\end{proof}

\subsection{The special cases $n = 2$ and $n = 3$}

The equation~\eqref{ThePbl} has no solution if \(n = 1\). Consider now the simplest cases \(n = 2\) and \(n = 3\).

\begin{lemma}\label{Lemn2}
\mbox{}

\begin{itemize}
\item[(i)] For \(n = 2\), a~ pair \((\overline{c_1},\overline{c_2})\) is a~ solution of equation~\eqref{ThePbl}
if and only if
\[
(\overline{c_1},\overline{c_2}) = (\overline{0},\overline{0}).
\]

\item[(ii)] For \(n = 3\), the triple \((\overline{c_1},\overline{c_2},\overline{c_3})\) is a~ solution of equation~\eqref{ThePbl}
if and only if
\[
(\overline{c_1},\overline{c_2},\overline{c_3}) = (\overline{1},\overline{1},\overline{1}).
\]
\end{itemize}
\end{lemma}

\begin{proof}
Part(i).
This follows from the equality
\[
\begin{pmatrix}
    \overline{c_{1}} & \overline{1} \\
    \overline{1} & \overline{0}
    \end{pmatrix}
\begin{pmatrix}
    \overline{c_{2}} & \overline{1} \\
    \overline{1} & \overline{0}
    \end{pmatrix}
=
\begin{pmatrix}
    \overline{c_{1}c_{2}+1} & \overline{c_{1}} \\
    \overline{c_{2}} & \overline{1}
    \end{pmatrix}.
\]
Part (ii).
This follows from the equality
\begin{align*}
\begin{pmatrix}
    \overline{c_{1}} & \overline{1} \\
    \overline{1} & \overline{0}
    \end{pmatrix}
\begin{pmatrix}
    \overline{c_{2}} & \overline{1} \\
    \overline{1} & \overline{0}
    \end{pmatrix}
\begin{pmatrix}
    \overline{c_{3}} & \overline{1} \\
    \overline{1} & \overline{0}
    \end{pmatrix}
&=
\begin{pmatrix}
    \overline{c_{1}c_{2}+1} & \overline{c_{1}} \\
    \overline{c_{2}} & \overline{1}
    \end{pmatrix}
\begin{pmatrix}
    \overline{c_{3}} & \overline{1} \\
    \overline{1} & \overline{0}
    \end{pmatrix}
\\
&=
\begin{pmatrix}
    \overline{c_{1}c_{2}c_{3}+c_{1}+c_{3}} & \overline{c_{1}c_{2}+1} \\
    \overline{c_{2}c_{3}+1} & \overline{c_{2}}
    \end{pmatrix}.
\end{align*}
Hence the result.
\end{proof}

\subsection{Inductive construction}
We now give an inductive procedure for the construction of all the solutions of equation~\eqref{ThePbl}
starting from the simplest case \(n = 2\) and the corresponding solution \((\overline{0},\overline{0})\).

\begin{proposition}\label{IndPro}
An \(n\)-tuple \((\overline{c_1},\ldots,\overline{c_n})\) is a~ solution of equation~\eqref{ThePbl}
if and only if
the sequence \((\overline{c_1},\ldots,\overline{c_n})\) can be obtained from
\((\overline{0},\overline{0})\) by applying the operations (a) and (b) in any order.
\end{proposition}

\begin{proof}
The ``if'' part follows from Lemma~\ref{FirstLem}.

Conversely, given a~ solution \((\overline{c_1},\ldots,\overline{c_n})\) of equation~\eqref{ThePbl},
one has the following two possibilities.
\begin{enumerate}
\item[A)] One has \(\overline{c_{i}} = \overline{0}\) for all \(i\).
Then, \(n\) is even and \((\overline{c_{1}},\ldots,\overline{c_{n}})\)
is obtained from \((\overline{0},\overline{0})\) by a~ sequence of \(\frac{n-2}{2}\) operations of type (b).

\item[B)] \(\overline{c_{i}} = \overline{1}\) for some \(i\), \(1 \leq i\leq n\).
Then, the inverse of the operation of type (a) centered at \(i\)
can be applied to \((\overline{c_{1}},\ldots,\overline{c_{n}})\).
This results in an \((n-1)\)-tuple
\((\overline{c_{1}},\ldots,\overline{c_{i-1}+1},\overline{c_{i+1}+1},\ldots,\overline{c_{n}})\).
The same computation as in the proof of Lemma~\ref{FirstLem} implies
that this \((n-1)\)-tuple
is a~ solution of equation~\eqref{ThePbl}.
We conclude by the induction assumption.
\qed
\end{enumerate}
\end{proof}

Let us mention that there exists an analog of Proposition~\ref{IndPro} in the case of non-negative integer solutions of~\eqref{CoCEqI}, see~\cite[Thm 3.1]{Cun}.
\subsection{End of the proof of Theorem~\ref{MaThm}}
We will need the following combinatorial interpretation of operations (a) and (b).
Let \((\overline{c_{1}},\ldots,\overline{c_{n}})\) be
a sequence corresponding to a~ \((3|4)\)-dissection of a~ convex \(n\)-gon,
then the result of either operation is again a~ sequence corresponding to a~ \((3|4)\)-dissection
of a~ convex \((n+1)\)-gon or \((n+2)\)-gon, respectively.
\begin{enumerate}
\item[(i)]
To a~ \((3|4)\)-dissection, operation (a) glues a~ triangle on the segment \((i,i + 1)\).
\item[(ii)]
Operation (b) glues a~ quadrilateral on the segment \((i,i + 1)\).
\end{enumerate}
We are ready to complete the proof of Theorem~\ref{MaThm}.

Part (i).
The induction basis consists of two cases, \(n = 3\) and \(n = 4\).
For \(n = 3\), the quiddity of a~ \((3|4)\)-dissection of any triangle is \((\overline{1},\overline{1},\overline{1})\) which is a~ solution of equation~\eqref{ThePbl}.
For \(n = 4\), the quiddity of a~ \((3|4)\)-dissection of any quadrilateral is \((\overline{1},\overline{0},\overline{1},\overline{0})\) (quadrilateral cut into two triangles)
or \((\overline{0},\overline{0},\overline{0},\overline{0})\) (quadrilateral alone) and it follows from Lemma~\ref{FirstLem} that they are solutions of equation~\eqref{ThePbl}.

Assume that an \(n\)-tuple \((\overline{c_1},\ldots,\overline{c_n})\) is the quiddity of a
\((3|4)\)-dissection of a~ convex \(n\)-gon.
Every \((3|4)\)-dissection has (at least) one exterior triangle (such a~ triangle is sometimes called an ``ear'' in the literature), or quadrilateral.
Cutting this exterior piece, one obtains either an \((n-1)\)-tuple or an \((n-2)\)-tuple
of elements of \(\Z/2\Z\)
which is the quiddity of a~ \((3|4)\)-dissection of a~ convex \((n-1)\)-gon or a~ convex \((n-2)\)-gon.
The result then follows from Lemma~\ref{FirstLem} and the induction assumption.

Part (ii).
For \(n = 3\), a~ triple \((\overline{c_1},\overline{c_2},\overline{c_3})\) is a~ solution of equation~\eqref{ThePbl}
if and only if \((\overline{c_1},\overline{c_2},\overline{c_3}) = (\overline{1},\overline{1},\overline{1})\),
which corresponds to a~ triangle.
Similarly to Lemma~\ref{Lemn2}, one shows the following: for \(n = 4\), the (cyclically ordered) solutions are \((\overline{1},\overline{0},\overline{1},\overline{0})\)
and \((\overline{0},\overline{0},\overline{0},\overline{0})\),
already considered in Example~\ref{TheEx}.

Assume that a~ sequence \((\overline{c_1},\ldots,\overline{c_n})\) is a~ solution
of equation~\eqref{ThePbl}, and let us show that it is the quiddity of a
\((3|4)\)-dissection of a~ convex \(n\)-gon.
By Proposition~\ref{IndPro}, this sequence is obtained from
\((\overline{0},\overline{0})\) by a~ sequence of the surgery operations (a) and (b).

If \(\overline{c_{i}} = \overline{1}\) for some \(i\),
where \(0 \leq i \leq n\), then, by the induction assumption, the sequence
\[
(\overline{c_{1}},\ldots,\overline{c_{i-1}+1},\overline{c_{i+1}+1},\ldots,\overline{c_{n}})
\]
is the quiddity of a~ \((3|4)\)-dissection
of a~ convex \((n-1)\)-gon.
Therefore,
\((\overline{c_1},\ldots,\overline{c_n})\) is the quiddity of a~ \((3|4)\)-dissection
of a~ convex \(n\)-gon, obtained from this \((3|4)\)-dissection by the gluing of a~ triangle.

If \(\overline{c_{i}} = \overline{c_{i+1}} = \overline{0}\), then
the sequence is of the form
\(
(\overline{c_{1}},\ldots,\overline{c_{i-1}},\overline{0},\overline{0},\overline{c_{i+2}},\ldots,\overline{c_{n}}).
\)
By the induction assumption,
\(
(\overline{c_{1}},\ldots,\overline{c_{i-1}},\overline{c_{i+2}},\ldots,\overline{c_{n}})
\)
 is the sequence associated to a~ \((3|4)\)-dissection
of a~ convex \((n-2)\)-gon.
Therefore
\((\overline{c_1},\ldots,\overline{c_n})\) is the quiddity of a~ \((3|4)\)-dissection
of a~ convex \(n\)-gon, obtained from this \((3|4)\)-dissection by the gluing of a~ quadrilateral.
Theorem~\ref{MaThm} is proved. \hfill \qed

\begin{remark}
Part (ii) of Theorem~\ref{MaThm} can be strengthened.
Let \((c_1,\ldots,c_n)\) be a~ solution of equation~\eqref{ThePbl}.
Assume that at least one element \(\overline{c_{i}}\) of \(\Z/2\Z\) is different from \(\overline{0}\)
(i.e., not all of \(c_i\) are even).
It turns out that, under this assumption, \((\overline{c_{1}},\ldots,\overline{c_{n}})\) is the quiddity of a~ triangulation of a~ convex \(n\)-gon.
For example, the following two \((3|4)\)-dissections have the same quiddity \((\overline{1},\overline{1},\overline{1},\overline{0},\overline{0})\).
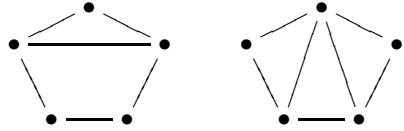
\[
\qquad
\xymatrix @!0 @R = 0.50cm @C = 0.5cm
{
&&\bullet\ar@{-}[rrd]\ar@{-}[lld]&
\\
\bullet\ar@{-}[rdd]\ar@{-}[rrrr]&&&& \bullet\ar@{-}[ldd]\\
\\
&\bullet\ar@{-}[rr]&& \bullet
}
\qquad
\xymatrix @!0 @R = 0.50cm @C = 0.5cm
{
&&\bullet\ar@{-}[rrd]\ar@{-}[lld]\ar@{-}[lddd]\ar@{-}[rddd]&
\\
\bullet\ar@{-}[rdd]&&&& \bullet\ar@{-}[ldd]\\
\\
&\bullet\ar@{-}[rr]&& \bullet
}
\]
The proof of this strengthened statement is very similar to that of Theorem~\ref{MaThm}, Part (ii).
It uses the following idea: if an \((n-1)\)-tuple, obtained by applying
the operation inverse to (a) centered at \(i\) to an \(n\)-tuple \((\overline{c_{1}},\ldots,\overline{c_{n}})\), contains only
 \(\overline{0}\), then the \((n-1)\)-tuple obtained by applying the operation inverse to (a)
 and centered at \(i+1\) to the \(n\)-tuple \((\overline{c_{1}},\ldots,\overline{c_{n}})\)
 contains an element different from \(\overline{0}\).
\end{remark}

\begin{remark}
Let us mention that Part (i) of Theorem~\ref{MaThm} can be deduced
from the combinatorial model and results of~\cite[Thm 7.3]{CH}.
One can also deduce from this model that
every quiddity \((\overline{c_1},\ldots,\overline{c_n})\) coming from a
\((3|4)\)-dissection can be lifted to an integer solution \(({c_1},\ldots,{c_n})\) of equation~\eqref{CoCEqI}.
\end{remark}
\medskip

\textbf{Concluding remark.}
Equation~\eqref{ThePbl} naturally extends to arbitrary
principal congruence subgroups \(\Gamma(N)\) in \(\SL(2,\Z)\),
and it would be interesting to find a~ combinatorial description of the set of solutions
in the general situation.

\medskip

\textbf{Acknowledgments.}
I am grateful to Valentin Ovsienko for stating the problem
and for his guidance during this work. I am also grateful to Michael Cuntz and Thorsten Holm for valuable comments.

\Editinfo{%
    November 17, 2019}{%
    December 20, 2019}{%
    Sophie Morier-Genoud
}

\end{paper}
\begin{references}

\refer{Paper}{Alp}
\Rtitle{Notes: PSL$_2$(Z) $=$ Z$_2 \ast$Z$_3$}
\Rauthor{R. Alperin}
\Rjournal{Amer. Math. Monthly}
\Rvolume{100}
\Rnumber{4}
\Rpages{385-386}
\Ryear{1993}

\refer{Paper}{ARS}
  \Rtitle{Friezes}
  \Rauthor{I. Assem, C. Reutenauer, D. Smith}
  \Rjournal{Adv. Math.}
  \Rvolume{225}
  \Rnumber{6}
  \Rpages{3134-3165}
  \Ryear{2010}

\refer{Paper}{BR}
  \Rtitle{SL$_k$-tilings of the plane}
  \Rauthor{F. Bergeron, C. Reutenauer}
  \Rjournal{Illinois J. Math.}
  \Rvolume{54}
  \Rnumber{1}
  \Rpages{263-300}
  \Ryear{2010}

\refer{Paper}{BHJ}
  \Rtitle{Generalized frieze pattern determinants and higher angulations of polygons}
  \Rauthor{C. Bessenrodt, T. Holm, P. J{\o}rgensen}
  \Rjournal{J. Combin. Theory Ser. A}
  \Rvolume{123}
  \Rpages{30-42}
  \Ryear{2014}

\refer{Paper}{CO}
  \Rtitle{Rotundus: triangulations, Chebyshev polynomials, and Pfaffians}
  \Rauthor{C. H. Conley, V. Ovsienko}
  \Rjournal{Math. Intelligencer}
  \Rvolume{40}
  \Rnumber{3}
  \Rpages{45-50}
  \Ryear{2018}

\refer{Paper}{CoCo}
  \Rtitle{Triangulated polygons and frieze patterns}
  \Rauthor{J. H. Conway, H. S. M. Coxeter}
  \Rjournal{Math. Gaz.}
  \Rvolume{57}
  \Rnumber{400}
  \Rpages{87-94}
  \Ryear{1973}

\refer{Paper}{CoCo2}
  \Rtitle{Triangulated polygons and frieze patterns}
  \Rauthor{J. H. Conway, H. S. M. Coxeter}
  \Rjournal{Math. Gaz.}
  \Rvolume{57}
  \Rnumber{401}
  \Rpages{175-183}
  \Ryear{1973}

\refer{Paper}{Cox}
  \Rtitle{Frieze patterns}
  \Rauthor{H. S. M. Coxeter}
  \Rjournal{Acta Arith.}
  \Rvolume{18}
  \Rpages{297-310}
  \Ryear{1971}

\refer{Paper}{Cun}
  \Rtitle{A combinatorial model for tame frieze patterns}
  \Rauthor{M. Cuntz}
  \Rjournal{Münster J. Math.}
  \Rvolume{12}
  \Rnumber{1}
  \Rpages{49-56}
  \Ryear{2019}

\refer{Paper}{CH}
  \Rtitle{Frieze patterns over integers and other subsets of the complex numbers}
  \Rauthor{M. Cuntz, T. Holm}
  \Rjournal{J. Comb. Algebra}
  \Rvolume{3}
  \Rnumber{2}
  \Rpages{153-188}
  \Ryear{2019}

\refer{Paper}{Hen}
  \Rtitle{Coxeter friezes and triangulations of polygons}
  \Rauthor{C.-S. Henry}
  \Rjournal{Amer. Math. Monthly}
  \Rvolume{120}
  \Rnumber{6}
  \Rpages{553-558}
  \Ryear{2013}

\refer{Paper}{HJ}
  \Rtitle{A $p$-angulated generalisation of Conway and Coxeter's theorem on frieze patterns}
  \Rauthor{T. Holm, P. J{\o}rgensen}
  \Rjournal{Int. Math. Res. Not. IMRN}
  \Rvolume{1}
  \Rpages{71-90}
  \Ryear{2020}

\refer{Paper}{Sop}
  \Rtitle{Coxeter's frieze patterns at the crossroads of algebra, geometry and combinatorics}
  \Rauthor{S. Morier-Genoud}
  \Rjournal{Bull. Lond. Math. Soc.}
  \Rvolume{47}
  \Rnumber{6}
  \Rpages{895-938}
  \Ryear{2015}

\refer{Paper}{Sop1}
  \Rtitle{Counting Coxeter's friezes over a~ finite field via moduli spaces}
  \Rauthor{S. Morier-Genoud}
  \Rjournal{Algebr. Comb.}
  \Rvolume{4}
  \Rnumber{2}
  \Rpages{225-240}
  \Ryear{2021}

\refer{Paper}{MO}
  \Rtitle{Farey boat: continued fractions and triangulations, modular group and polygon dissections}
  \Rauthor{S. Morier-Genoud, V. Ovsienko}
  \Rjournal{Jahresber. Dtsch. Math.-Ver.}
  \Rvolume{121}
  \Rnumber{2}
  \Rpages{91-136}
  \Ryear{2019}

\refer{Other}{OEIS}
  \Rtitle{OEIS Foundation Inc., The On-Line Encyclopedia of Integer Sequences}
  \Rnote{http://oeis.org}

\refer{Paper}{O}
  \Rtitle{Partitions of unity in SL(2, Z), negative continued fractions, and dissections of polygons}
  \Rauthor{V. Ovsienko}
  \Rjournal{Res. Math. Sci.}
  \Rvolume{5}
  \Rnumber{2}
  \Rpages{paper no. 21, 25 pp}
  \Ryear{2018}

\refer{Book}{R}
\Rauthor{R. Rankin}
\Rtitle{The modular group and its subgroups}
\Rpublisher{The Ramanujan Institute, Madras}
\Ryear{1969}

\refer{Paper}{Zag}
\Rauthor{D. Zagier}
  \Rtitle{Nombres de classes et fractions continues}
  \Rjournal{Journées Arithmétiques de Bordeaux. Astérisque}
  \Rvolume{24-25}
  \Rnumber{0}
  \Rpages{81-97}
  \Ryear{1975}

\end{references}
